\documentclass[a4paper,11pt]{amsart}
\usepackage{amssymb}
\usepackage{amsmath}
\usepackage{graphicx}

\font\eufm=eufm10
\def\frak#1{\hbox{\eufm#1}}
\newcommand{\bd}{\begin{document}}
\newcommand{\ed}{\end{document}}
\newcommand{\be}{\begin{enumerate}}
\newcommand{\ee}{\end{enumerate}}
\newcommand{\bi}{\begin{itemize}}
\newcommand{\ei}{\end{itemize}}
\newcommand{\ba}{\begin{array}}
\newcommand{\ea}{\end{array}}
\newcommand{\vs}{\vspace*{0.3\baselineskip}}
\newcommand{\vsm}{\vspace*{-0.3\baselineskip}}
\newcommand{\kom}[1]{{\em #1}\newline}
\newtheorem{defi}{Definition}[section]
\newtheorem{tw}[defi]{Theorem}
\newtheorem{prop}[defi]{Proposition}
\newtheorem{lem}[defi]{Lemma}
\newtheorem{re}[defi]{Remark}
\newtheorem{col}[defi]{Corollary}
\newtheorem{ex}[defi]{Examples}
\newtheorem{ex1}[defi]{Example}
\newtheorem{zad}{Exercise}[section]
\newtheorem{zal}{Assumptions}[section]
\newtheorem{assumpt}[defi]{Assumptions}
\newcommand{\Om}{\Omega}
\newcommand{\om}{\omega}
\newcommand{\G}{\Gamma}
\newcommand{\D}{\Delta}
\renewcommand{\d}{\delta}
\newcommand{\ga}{\gamma}
\newcommand{\eps}{\epsilon}
\newcommand{\ove}{\overline}
\newcommand{\ms}{\oplus}
\newcommand{\mt}{\otimes}
\newcommand{\dz}{\wedge}
\newcommand{\lra}{\longrightarrow}
\newcommand{\sign}{\mbox{$ sgn $}}
\newcommand{\rra}{\rightrightarrows}
\newcommand{\rel}{\mbox{$\,$\rule[0.5ex]{1.1em}{0.2pt}$\triangleright\,$}}
\newcommand{\dow}{\hspace*{\fill}\rule{1.6ex}{1.6ex}\hspace*{1em}}
\newcommand{\dowl}{\hspace*{\fill}\rule{1ex}{1ex}\hspace*{1em}}
\newcommand{\sd}{\hspace{0.3ex}\tiny{\rhd\mbox{\hspace{-2ex}}<}\hspace{0.3ex}}
\newcommand{\mmt}[2]{\mbox{$\vphantom{}_{#1}\times_{#2}$}}
\newcommand{\gotg}{\frak g}
\newcommand{\gota}{\frak a}
\newcommand{\gotb}{\frak b}
\newcommand{\gotc}{\frak c}
\newcommand{\gothh}{\frak h}
\newcommand{\gott}{\frak t}
\newcommand{\hd}{\hat{\d}}
\newcommand{\oml}{\Omega_L^{1/2}}
\newcommand{\omr}{\Omega_R^{1/2}}
\newcommand{\omh}{\Omega^{1/2}}
\newcommand{\lo}{\lambda_0}
\newcommand{\ro}{\rho_0}
\newcommand{\lma}{\Lambda^{max}}
\newcommand{\timh}{\times_h}
\newcommand{\Gd}{\G^{(2)}}
\newcommand{\el}{e_L}
\newcommand{\er}{e_R}
\newcommand{\GG}{\G_1\times\G_2}
\newcommand{\gdot}{\hspace{-0.1em}\cdot\hspace{-0.1em}}
\newcommand{\tran}{\frown\hspace{-2.2ex}|\hspace{1.9ex}}
\newcommand{\la}[2]{\Lambda_{#1#2}}
\newcommand{\kad}{ad^{\#}}
\newcommand{\wl}[1]{\vphantom{X}_{#1}{\G}}
\newcommand{\te}{\tilde{e}}
\newcommand{\notka}[1]{}
\newcommand{\sA}{\mbox{$\mathcal A$}}
\newcommand{\sT}{\mbox{$\mathcal T$}}
\newcommand{\sB}{\mbox{$\mathcal B$}}
\newcommand{\sF}{\mbox{$\mathcal F$}}
\newcommand{\sO}{\mbox{$\mathcal O$}}
\newcommand{\sD}{\mbox{$\mathcal D$}}
\newcommand{\sS}{\mbox{$\mathcal S$}}
\newcommand{\sY}{\mbox{$\mathcal Y$}}
\newcommand{\hY}{\mbox{$\hat{Y}$}}
\newcommand{\hS}{\mbox{$\hat{S}$}}
\newcommand{\hX}{\mbox{$\hat{X}$}}
\newcommand{\dif}{differential }
\newcommand{\gru}{groupoid }
\newcommand{\grus}{groupoids }
\newcommand{\ti}{\tilde}
\newcommand{\halden}{half density }
\newcommand{\haldens}{half densities }
\renewcommand{\top}{topological }
\newcommand{\Setrel}{\mbox{\rm SetRel}}
\newcommand{\cstardwa}{\mbox{$C^*_r(\Gamma\times\Gamma)$}}
\begin{document}
\title{Short and biased introduction to groupoids}
\author{Piotr Stachura}
\address{Faculty of Applied Informatics and Mathematics, Warsaw University of Life Sciences-SGGW,
ul Nowoursynowska 166, 02-787 Warszawa, Poland,  
e-mail: stachura@fuw.edu.pl, piotr\_stachura1@sggw.pl}
\date{}
\begin{abstract} The algebraic part  of approach to groupoids started by S. Zakrzewski is presented.
\end{abstract}

\maketitle
\section{Introduction}
These notes  present  the approach to groupoids started by S. Zakrzewski \cite{SZ}.
They are  focused on {\em algebraic structure} and any topological or differential topics are completely ignored 
(see \cite{SZ,PS,BunPs,CatDheWei}). I hope they might be useful for people who want to be quickly introduced to the subject.

There are some new observations  that includes formal properties of a category of groupoids considered here and a new 
(but  equivalent to the standard one) definition of  a groupoid action.

The next section contains definitions (the standard one and given by S. Zakrzewski) and main examples of groupoids; it also  fixes notation.
In the third one, Zakrzewski morphisms are defined and examples are given; some formal properties 
of the resulting category of groupoids are investigated; finally, in the last section we discuss relations of  Zakrzewski 
morphisms to groupoid actions -- they give alternative view on Zakrzewski morphisms.

The bibligraphy includes only articles  directly related to the presented approach. An interested reader may consult e.g. \cite{McK, McK1}
for a more complete list. 
\section{Groupoids: definitions and examples}
Groupoids are special categories and the shortest definition is:
\begin{defi}\label{grup-cat}
A groupoid is a small category with invertible morphisms.
\end{defi}

Thus we have a set  of objects $E$ {\em (small category!)} and a set of morphisms $\G$. We will identify an object $e\in E$ with 
the corresponding identity morphism ${\rm id}_e\in\G$. After this identification $E$ is a subset of $\G$.
Each morphism 
$\gamma\in\G$ has its {\em domain (source)}: $d(\gamma)\in E$ and its {\em range (target)}: $r(\gamma)\in E$. 
Composition $\gamma_1\gamma_2$ of two morphisms is possible iff $d(\gamma_1)=r(\gamma_2)$; then
$d(\gamma_1 \gamma_2)=d(\gamma_2)$ and $r(\gamma_1 \gamma_1)=r(\gamma_1)$. Every 
morphism $\gamma$ has its {\em inverse} $s(\gamma)$ with $s(s(\gamma))=\gamma$ and $d(s(\gamma))=r(\gamma)$. 
Composition is associative $\gamma_1(\gamma_2\gamma_3)=(\gamma_1\gamma_2)\gamma_3$ (in a sense that if one side is defined the 
other also and then they are equal) and inverse reverses the order of composition: $s(\gamma_1 \gamma_2)=s(\gamma_2) s(\gamma_1)$. 

Groupoid $\G$ with the set of identities $E$ will be denoted by $\G\rightrightarrows E$; the set of identities will be also 
denoted by $\G^{(0)}$ and the set {\em of composable pairs} by 
$\Gamma^{(2)}:=\{(\gamma_1,\gamma_2) \in \Gamma\times \G: d(\gamma_1)=r(\gamma_2)\}$.
For $e\in E$ the set $d^{-1}(e)\cap r^{-1}(e)=:\Gamma_e$
is a group called the {\em isotropy group} of $e$. 
Union of all isotropy groups is called {\em the isotropy group bundle} (of $\G$).
A relation on $E$ defined by:
$e_1\sim e_2 \iff \exists \gamma\,:\,e_1=r(\gamma), e_2=d(\gamma)$ is an equivalence relation and its classes 
are {\em orbits} of $\G$; the set of orbits will be denoted by $E/\Gamma$. 
If $\sO\in E/\Gamma$ is an orbit then the set $d^{-1}(\sO)$ ($=r^{-1}(\sO)$) is called a {\em transitive component} of $\G$. Every groupoid is 
a (disjoint) union of its transitive components. 
{\em Subgroupoid} of a groupoid $\G$ will mean a subset $G\subset\G$  that satisfies two conditions:
 $g\in G\Rightarrow s(g)\in G$ and $(g_1,g_2)\in\G^{(2)}\cap G\Rightarrow g_1 g_2\in G$; subgroupoid $G\subset \G$ is {\em wide} iff $E\subset G$.

\subsection{Examples}

{\em Group bundles and groups.} If $d=r$ then each element $e\in E$ is an orbit and each set $d^{-1}(e)$ is a transitive component and a group; 
thus $\G$ is a disjoint union of groups. In particular, if $E$ consists 
of only one element then $\G$ is a group.

{\em Sets.} These are special group bundles: every group is trivial. In this situation $d=r=s=id$,  $\Gamma^{(2)}=diag(\Gamma\times\Gamma)$, 
$\G^{(0)}=\G$ and $\gamma\gamma=\gamma$ for every $\gamma\in\G$.

{\em Pair groupoids.} These are transitive groupoids with trivial isotropy groups (see lemma \ref{tran-str}):  
$\G=X\times X \,,\, E:=diag(X\times X)\simeq X$; $d(x,y):=(y,y)\,,\,r(x,y):=(x,x)$. Composition is defined as $(x,y)(y,z):=(x,z)$ and 
inverse $s(x,y):=(y,x)$. For a set $X$ this groupoid will be denoted by $X^2$.

{\em Equivalence relations.} If $R\subset X\times X$ is an equivalence relation, 
then it is a groupoid with operations inherited from the pair groupoid $X^2$. Its orbits are equivalence classes and 
transitive components are pair groupoids based on equivalence classes.

{\em Products of groups and pair groupoids.} Let $G$ be a group and $X$ a set. Define $\G:=X\times G\times X$, 
$E:=\{(x,e,x): x\in X\}$ (here $e$ is the neutral element of $G$)  and  $d(x,g,y):=(y,e,y)$, $r(x,g,y):=(x,e,x)$. 
Multiplication is defined as $(x,g,y)(y,h,z):=(x,gh,z) $ and inverse $s(x,g,z):=(z,g^{-1},x)$. 
(These are really cartesian products of groups and pair groupoids.)

{\em Transformation groupoids}. Let a group $G$ acts on a set $X$, i.e. there is a mapping
 $G\times X\ni (g,x)\mapsto gx\in X$   that satisfies usual conditions. Then the transformation groupoid is 
$\Gamma:=G\times X$ with the set of units  $E:=\{(e,x): x\in X\}\simeq X$; domain and range maps: 
$d(g,x):=(e,x)$, $r(g,x):=(e,gx)$; inverse: $s(g,x):=(g^{-1}, gx)$ and composition: $(g,hx) (h, x):=(gh,x)$.
Orbits of $\G$ are orbits of the action and transitive components are groupoids from the last example with groups being
stabilizers and pair groupoids are based on orbits.

{\em From purely algebraic}  point of view (lemma \ref{tran-str}) every groupoid can be built from groups and pair 
groupoids by products and disjoint union.

\subsection{Alternative definition.}

An alternative definition of groupoids was given by S. Zakrzewski in \cite{SZ}. It is based on relations so let us now recall some
 basic facts and introduce notation. 

Let $X, Y$ be sets. A {\em relation $r$ from $X$ to $Y$} is a triple $(X,Y;Gr(r))$, where $Gr(r)\subset Y\times X$.
A relation from $X$ to $Y$ will be denoted by $r: X\rel Y$ {\em (note the type of arrow!)}. We will usually identify a relation 
$r:X\rel Y$ with its graph and will write $(y,x)\in r$ instead of $(y,x)\in Gr(r)$ since it shouldn't cause any confusion.
For a relation $r:X\rel Y$ its {\em transposition} $r^T:Y\rel X$ is defined by $(x,y)\in r^T\iff (y,x)\in r$; 
the {\em domain} of $r$ is the set $D(r):=\{x\in X: \exists y\in Y \,(y,x)\in Gr(r)\}$ and the {\em image} of $r$ is the  set 
$Im(r):=\{y\in Y: \exists x\in X \, (y,x)\in Gr(r)\}$
A {\em composition} of  relations $r:X\rel Y$ and  $s:Y\rel Z$ is the relation  $sr:X\rel Z$ defined by
$Gr(sr):=\{(z,x): \exists \, y\in Y: (z,y)\in s\,,\,(y,x)\in r\}$. Sets with relations as morphisms form a category. 
{\em Cartesian product} of relations is defined in a natural way:
$r:X\rel Y$ and $r_1:X_1\rel Y_1$ then $r\times r_1: X\times X_1\rel Y\times Y_1$ is equal to:
$r\times r_1:=\{(y,y_1,x,x_1): (y,x)\in r\,,\,(y_1,x_1)\in r_1\}$ ({\em but this is not a product in a categorical sense}).
In the following the symbol
$\{1\}$ denotes  one point set and $\sigma:X\times Y\ni(x,y)\mapsto (y,x)\in Y\times X$ a flip.
\begin{defi}\label{def-grup}\cite{SZ} Groupoid  $\Gamma\rightrightarrows E$ is a quadrupole  $(\Gamma,m,s,e)$, where  $\G$ is a set, 
$m: \G\times \G\rel \G$, 
$e:\{1\}\rel \Gamma$ and $s: \G\rel \G$ are relations satisfying the following conditions:
$$m(m\times id)=m(id\times m)\,\,\,,\,\,m(e\times id)=m(id\times e)=id\vspace{-1ex}$$
$$s^2=id\,\,,\,\,s m=m \sigma  (s\times s)\vspace{-1ex}$$
$$\forall \gamma\in\G \,\,\emptyset\neq m(s(\gamma),\gamma)\subset E:=Im(e)\vspace{-1ex}$$
\end{defi}
It follows that this  definition is equivalent to def. \ref{grup-cat}, in particular $m$ is in fact {\em a mapping} defined on 
$\Gamma^{(2)}\subset\Gamma\times\Gamma$; see \cite{SZ} for proofs and details. 
We will sometimes use $D(m)$ for $\G^{(2)}$ and 
{\em from now on we will use different notation for target and source mappings:}\vspace{-1ex}
$$e_L(\gamma):=m(\gamma,s(\gamma))=r(\gamma)\,\,,\,\,\,e_R(\gamma):=m(s(\gamma),\gamma)=d(\gamma)$$

\subsection{Operations on groupoids}

{\em Restriction.} If $\Gamma$ is a groupoid and $F\subset E$ then the set 
$e_L^{-1}(F)\cap e_R^{-1}(F)$ is a groupoid with $F$ as a set of identities  -- the restriction of $\Gamma$ to $F$.

{\em Disjoint union.} Disjoint union of groupoids is a groupoid. Any groupoid is a disjoint union 
of its transitive components and (essentially) this is the only way for a groupoid to be a disjoint union:
\begin{lem}\label{union} Let $\G$ be a groupoid and $\{\G_\lambda\,,\,\lambda\in \Lambda\}$ a family of subgroupoids of $\G$.  If
$\displaystyle \Gamma=\bigcup_{\lambda\in\Lambda}\G_\lambda$ and $\Gamma_\lambda\cap\Gamma_{\lambda'}=\emptyset$ for $\lambda\neq\lambda'$ then 
each $\Gamma_\lambda$ is a union of transitive components (of $\Gamma$).\dowl
\end{lem}

{\em Cartesian product.} If $\Gamma_1$ and $\Gamma_2$ are groupoids then $\Gamma_1\times \Gamma_2$ 
is a groupoid. The set of units is $E_1\times E_2$ and the multiplication relation is 
$(m_1\times m_2)(id_{\Gamma_1}\times \sigma\times id_{\Gamma_2})$, where $\sigma(\gamma_2,\gamma_1)=(\gamma_1,\gamma_2)$ is a flip. 

The following proposition describes a structure of transitive groupoids; although morphisms of groupoids have not been defined yet, the 
{\em `` isomorphism''} used in the lemma means  
{\em a  bijective mapping preserving all components of the groupoid structure.}
\begin{prop}\label{tran-str}Let $\Gamma$ be a transitive groupoid. Let us choose $e_0\in E$ and let  $G$ be the isotropy 
subgroup of $e_0$. Let $p:E\rightarrow e_L^{-1}(e_0)$ be a section of $e_R$ such that
$p(e_0)=e_0$. The mapping:
$$E\times G\times E\ni(x,g,y)\mapsto s(p(x)) g p(y)\in\Gamma$$
is an isomorphism (groupoid structure on $E\times G\times E$ is the product one).  \dow
\end{prop}

\section{Morphisms of groupoids.}

If one considers groupoid as a special type of category,
the natural notion of morphism is a functor, i.e.  a map
$\Phi: \Gamma_1\rightarrow\Gamma_2$ that satisfies:
$\Phi(E_1)\subset E_2$; $\Phi(s(\gamma))=s(\Phi(\gamma))$ and, 
for any composable $(\gamma,\tilde{\gamma})\subset \Gamma^{(2)}_1$, 
$\Phi(\gamma \tilde{\gamma})=\Phi(\gamma)\Phi(\tilde{\gamma})$.
But the definition \ref{def-grup} suggests another possibility \cite{SZ}:
\begin{defi}\label{def-mor}{\em (Zakrzewski  morphisms.)}
Let $\Gamma\rightrightarrows E,\Gamma_1\rightrightarrows E_1$ be groupoids. A morphism from $\Gamma$ to $\Gamma_1$  is a relation
$h:\Gamma\rel\Gamma_1$ that satisfies:
$$ hm=m_1 (h\times h)\,\,,\,h s=s_1 h\,\,,\,h e=e_1$$
\end{defi}
It follows from this  definition  (see \cite{SZ}) that a morphism  $h:\G\rel \G_1$ determines:\\
a map  $\rho_h:E_1\rightarrow E$ called {\em the base map of $h$}
and for  each $e\in E_1$ two mappings of fibers:
$$ h_e^L: e_L^{-1}(\rho_h(e))\rightarrow e_L^{-1}(e)\,\,,\,\,h_e^R: e_R^{-1}(\rho_h(e))\rightarrow e_R^{-1}(e)$$
In particular the image of $\rho_h$ is a union of orbits, the domain of $h$ is a union of transitive components and 
the image of $h$ is a wide subgroupoid of $\Gamma_1$. If $\G',\G''\subset\G$ are transitive components their 
images by $h$ are disjoint: $h(\G')\cap h(\G'')=\emptyset$. 

As the name suggests, {\em groupoids with Zakrzewski morphisms form a category} i.e. 
composition of Zakrzewski  morphisms (as relations) is a Zakrzewski morphism.
In the following we will almost exclusively consider Zakrzewski  morphisms  
and they will be called  morphisms. The standard ones will be called functors. 

\begin{re}It should be pointed out that  Zakrzewski  morphisms are  not generalization of functors. 
A functor $\Phi: \Gamma\rightarrow\Delta$ is a Zakrzewski  morphism iff $\Phi$ restricted to 
$\Gamma^{(0)}$ is a bijection (onto $\Delta^{(0)}$).
\end{re}

\noindent
The image of a subgroupoid by a morphism is a subgroupoid and one can restrict morphisms to subgroupoids.
\begin{prop}
Let $h:\Gamma\rel\Delta$ be a morphism of groupoids and  $G,G_1\subset \G$ subgroupoids.
\begin{enumerate}
\item $h(G)\subset \Delta$ is a subgroupoid;
\item If $G\cap G_1=\emptyset$ then $h(G)\cap h(G_1)=\emptyset$;
\item $h|_G:G\rel h(G)$ is a morphism.
\item If $h$ is surjective and $G$ is a transitive component then $h(G)$ is a union of transitive components.
\end{enumerate}
\end{prop}
{\em Proof}: 1) is immediate from def. \ref{def-mor}; 2) If $(\delta,g)\in h$ and $(\delta, g_1)\in h$ for $g\in G\,,\,g_1\in G_1$ then
$e_L(g)=e_L(g_1)\in G\cap G_1$ -- impossible; 3) -- direct verification; 4) If $h$ is surjective then $\Delta=\bigcup h(\Gamma_\alpha)$, where
$\Gamma_\alpha$'s are transitive components of $\Gamma$. Now use point 2) and lemma \ref{union}.\\
\dow

\noindent
A morphism is determined by its value on any fiber in every transitive component (contained in its domain).
\begin{lem} Let $\Gamma\rightrightarrows E$  and $\Delta\rightrightarrows F$ be groupoids and  
$h,k:\Gamma\rel\Delta$ morphisms. Assume $\Gamma$ is transitive and for some $e\in E$: 
$h|_{\er^{-1}(e)}=k|_{\er^{-1}(e)}$. Then $h=k$.
\end{lem}
{\em Proof:} Let $(\delta,\gamma)\in h$; there exists $\gamma_1$ such that $\el(\gamma_1)=\el(\gamma)$ and $\er(\gamma_1)=e$ 
(transitivity of $\Gamma$). Let $\delta_1$ be the unique element such that $(\delta_1,\gamma_1)\in h $ and $\el(\delta_1)=\el(\delta)$. 
Let $\delta_2:= s(\delta) \delta_1$, then $(\delta_2,s(\gamma) \gamma_1)\in h$. Denoting $\gamma_2:=s(\gamma) \gamma_1$ we have:
$$(\delta_2,\gamma_2)\in h\,,\,(\delta_1,\gamma_1)\in h\,,\,(\delta,\gamma)\in h\,,\,\delta_1= \delta \delta_2\,,\,\gamma_1=\gamma\gamma_2.$$
By assumption $(\delta_1,\gamma_1)\in k$ and $(\delta_2,\gamma_2)\in k$ so $(s(\delta_2),s(\gamma_2))\in k$ and 
$(\delta;\gamma_1,s(\gamma_2))\in m'(k\times k)=k m$ ($m'$ is the multiplication in $\Delta$).
 In this way we get $(\delta,\gamma)\in k$. The opposite inclusion can be shown in  the same way.\vspace{-1ex}\\\dowl

\subsection{Examples of morphisms}
{\em Groups.}  If $\G_1,\G_2$ are groups, morphisms from $\G_1$ to $\G_2$ (as well as functors) are group homomorphisms.

{\em Sets.} Here functors from $\Gamma_1$ to $\Gamma_2$ are just mappings. 
Zakrzewski morphisms are also mappings but in a opposite direction, any morphism $h:\Gamma_1\rel\Gamma_2$ is equal to 
 $h=f^T=\{(x,f(x)) :x\in\Gamma_2\}$ for a map $f:\Gamma_2\rightarrow\Gamma_1$.

{\em  Regular representation (left). } For any groupoid the relation $l:\G\rel \G^2$ (pair groupoid)
defined by: \begin{equation}\label{left}l:=\{(g,h;k): (g;k, h)\in m\}\end{equation}
is a morphism.

{\em Transitive components.} If $\G_1\subset\G$ is a union of transitive components and $i: \G_1\rightarrow\G$ is the inclusion map,
then $i^T:\Gamma\rel \G_1$ is a morphism.

{\em Restriction of morphism to its domain.} If $h:\Gamma_1\rel\G_2$ is a morphism with a domain $D(h)$, 
then the relation $h|_{D(h)}: D(h)\rel\Gamma_2$  is a morphism.

{\em Wide subgroupoids.} If $\Gamma_1\subset \G$ is a wide subgroupoid (i.e. $E\subset\G_1$), the 
inclusion $i:\Gamma_1\rightarrow\Gamma$ is a morphism.

{\em Isotropy group bundle. } This is special type of the previous example.
Let $\G$ be a groupoid, $\G':=\bigcup_{e\in E}e_L^{-1}(e)\cap e_R^{-1}(e)$ its isotropy group bundle 
and $i:\Gamma'\rightarrow\G$ the inclusion. Then $i:\Gamma'\rel \Gamma$  is a morphism.

{\em For any groupoid}  a map: $\Gamma\ni \gamma\mapsto (e_L(\gamma),\er(\gamma))\in E\times E$ is a morphism to $E^2$; 
{\em this map is also  a  functor}.

{\em Cartesian products. } Let $\Gamma=\Gamma_1\times\Gamma_2$. Relations $i_1:=\{(g_1,e_2;g_1): g_1\in\G_1, e_2\in E_2\}$ and 
$i_2:=\{(e_1,g_2;g_2): g_2\in\G_2, e_1\in E_1\}$ are morphisms from $\G_1\rel \G_1\times \G_2$ and $\G_2\rel \G_1\times \G_2$. 
{\em But projections $\pi_1(\pi_2):\G_1\times\G_2\rightarrow \G_1(\G_2)$ are not morphisms.} 
So cartesian product of groupoids is not a product in categorical sense (it is rather like tensor product).

{\em Group actions. } If a group $G$ acts on a set $X$ (left action), then the action defines 
a morphism $h:G\rel X^2$:
$h=\{(g x,x;g) : g\in G\,,\,x\in X\}$
and, conversely, any morphism $G\rel X^2$ comes from an action. Unless $X$ is a one point set, this is not a mapping.

{\em Morphisms into groups.} Let $G$ be a group and $\Gamma\rra E$ a groupoid. Morphisms $h: \Gamma\rel G$ are given by a one point orbit 
$\{e_0\}\subset E$ and a homomorphism from the isotropy group of $e_0$ into $G$.

{\em Morphisms from groups.} Let $G$ be a group and $\Gamma\rra E$ a groupoid. Morphisms $h: G\rel \Gamma$ are just group homomorphisms 
to a group of bisections of $\G$ (see below).

\subsection{ Bisections and morphisms}

With any groupoid $\Gamma\rightrightarrows E$ there comes a group: {\em the group of bisections of $\Gamma$}.

\begin{defi}
A set $B\subset \Gamma$ is a bisection iff it is a section of left and right projection over $E$.
\end{defi}
\noindent For groups bisections are elements, for pair groupoids bisections are graphs of bijections.

Subsets of a groupoid can be ``multiplied'': if $A,B\subset \Gamma$ we define $AB:=\{m(a,b): a\in A, b\in B, (a,b)\in D(m)\}$.
This operation turns the set of bisections into  a group: neutral element is the set of identities and $B^{-1}=s(B)$.
This ``multiplication'' of subsets can be used to characterize bisections:
\begin{lem}
Let $\Gamma\rra E$ be a groupoid and $A\subset \Gamma$.
\begin{itemize} 
\item $A$ is a section of $e_R$ over $e_R(A)$ iff $As(A)\subset E$;
\item $A$ is a section of $e_L$ over $e_L(A)$ iff $s(A) A\subset E$;
\item $A$ is a bisection iff $s(A)A=As(A)=E$.
\end{itemize}
\vspace{-2ex}\dowl
\end{lem}
Bisections act on a groupoid by $\Gamma\ni\gamma\mapsto B\gamma:=\gamma'\gamma$, 
where $\gamma'$ is a unique element in $B$ with $\er(\gamma')=\el(\gamma)$ (i.e.
$\{B\gamma\}=B\{\gamma\}$ using ``multiplication'' of subsets). 
This action preserves right fibers i.e. $\er(B\gamma)=\er(\gamma)$ and maps left fibers into left fibers.
Morphisms act on bisections:
\begin{lem} \label{mor-and-bis} 
Let $h:\Gamma\rel \Gamma_1$ be a morphism and $B,B'\subset \Gamma$ bisections.  The set $h(B)$ is a bisection, $h(s(B))=s_1(h(B))$ and
$h(B B')=h(B) h(B')$ (in other words, $h$ defines a homomorphism of groups of bisections).
\dowl
\end{lem}
\noindent
Bisections act on a groupoid by ``adjoint'' action, this action is a morphism:
\begin{lem} Let $B,C\subset\Gamma$ be  bisections.
The mapping $Ad_B: \Gamma\ni\gamma\mapsto B{\gamma}s(B)\in\Gamma$ is a morphism; $Ad_{B C}=Ad_B Ad_C$. 
If $h: \Gamma\rel \Gamma_1$ is a morphism then $h\cdot Ad_B=Ad_{h(B)}\cdot h$.\dowl
\end{lem}
\subsection{Monomorphisms}
Now we give a complete characterization of monomorphisms, it turns out it is similar to group case. Let us begin with the following:
\begin{defi} The kernel of a morphism $h:\Gamma\rel \Delta$ is a set 
$\ker(h):=\{\gamma\in D(h): (\delta,\gamma)\in h\Rightarrow \delta\in \Delta^{(0)}\}$.
\end{defi}
The next lemma states that $\ker(h)$ is a subgroupoid contained in the isotropy group bundle of $\Gamma$ and its intersection 
with every isotropy group  over $\ker(h)\cap \Gamma^{(0)}$ is a normal subgroup.
\begin{lem} \label{ker-prop}Let $h:\Gamma\rel \Delta$ be a morphism. Then
\begin{enumerate}
\item $D(h)\cap \Gamma^{(0)}\subset \ker(h)$;
\item $\gamma\in \ker(h)\Rightarrow \el(\gamma)=\er(\gamma)$;
\item $\gamma\in\ker(h)\Rightarrow s(\gamma)\in\ker(h)$;
\item if $\gamma_1,\gamma_2\in \ker(h)$ and $(\gamma_1,\gamma_2)\in D(m)$ then $\gamma_1\gamma_2\in \ker(h)$;
\item if $\gamma\in \ker(h)$ and $(\gamma_1,\gamma)\in D(m)$ then $\gamma_1 \gamma s(\gamma_1)\in \ker(h)$;
\item if $B\subset \G$ is a bisection then $Ad_B (\ker(h))=\ker(h)$.
\end{enumerate}
\end{lem}

\noindent {\em Proof:} 1),  2) and 3) are obvious.

4) It is clear that $\gamma_1\gamma_2\in D(h)$. Let $(\delta,\gamma_1\gamma_2)\in h$ i.e. $(\delta;\gamma_1,\gamma_2)\in hm=m'(h\times h)$;
so exist $\delta_1,\delta_2$ such that $(\delta;\delta_1,\delta_2)\in m'$ and $(\delta_1,\gamma_1)\in h$, and $(\delta_2,\gamma_2)\in h$. Since 
$\gamma_1,\gamma_2\in \ker(h)$ it follows that $\delta_1=\delta_2=\delta\in\Delta^{(0)}$.

5) Let $m', s',etc$ be operations in $\Delta$. 
If $\gamma\in \ker(h)$ and $(\gamma_1,\gamma)\in D(m)$ then $\gamma_1,s(\gamma_1), \gamma_1 \gamma s(\gamma_1) \in D(h)$. 
Let $(\delta, \gamma_1 \gamma s(\gamma_1))\in h$;
there exists $\delta_1$ such that $\el'(\delta_1)=\el'(\delta)$ and $(\delta_1,\gamma_1)\in h$, so also $(s'(\delta_1),s(\gamma_1))\in h$. 
Therefore $(s'(\delta_1) \delta; s(\gamma_1), \gamma_1 \gamma s(\gamma_1))\in m'(h\times h)=h m$ and
so $(s'(\delta_1) \delta, \gamma s(\gamma_1))\in h$.
But then $h\ni (\el'(s'(\delta_1) \delta), \el(\gamma s(\gamma_1)))=(\er'(\delta_1),\el(\gamma))$ and, since $\gamma\in \ker(h)$, 
$(\er'(\delta_1),\gamma)\in h$, so also $(\er'(\delta_1),s(\gamma))\in h$. Therefore
$(s'(\delta_1) \delta;s(\gamma), \gamma s(\gamma_1))\in m'(h\times h)=hm$ and $(s'(\delta_1) \delta;s(\gamma_1))\in h$.
In this way we get $(\delta_1,\gamma_1)\in h$ and $(s'(\delta_1) \delta;s(\gamma_1))\in h$ but then  
$(\delta;\gamma_1,s(\gamma_1))\in m'(h\times h)=hm$ i.e. $(\delta,\el(\gamma_1))\in h$ and $\delta\in\Delta^{(0)}$.\\
6) is a direct consequence of 5).
\dowl
\begin{prop} A morphism $\phi:\Gamma\rel \Delta$ is monomorphism iff $\ker(\phi)=\Gamma^{(0)}$.
\end{prop}
{\em Proof:} $\Leftarrow\,\,$ 
Suppose $\phi$ is not a monomorphism. If $D(\phi)\neq \G$ then of course
$\ker(\phi)\neq \G^{(0)}$ so we can assume  that $D(\phi)=\G$ i.e. $\rho_\phi$ is surjective.
There exists a groupoid $K$ and morphisms $\psi_1,\psi_2:K\rel \Gamma$ such that $\psi_1\neq \psi_2$ and $\phi \psi_1=\phi \psi_2$.
The last equality means that $\rho_{\psi_1}\rho_\phi=\rho_{\psi_2}\rho_\phi$ so $\rho_{\psi_1}=\rho_{\psi_2}$,
in particular $D(\psi_1)=D(\psi_2)$.
Let $(\gamma_0,k_0)\in \psi_1\setminus\psi_2$.
Let $\er(\gamma_0)=:e_0$. Since $\phi\psi_1=\phi \psi_2$ for  every $\delta\in \phi(\gamma_0)$ there exists $g_\delta\neq \gamma_0$ such that 
$(\delta,\gamma_\delta)\in\phi$ and $(\gamma_\delta,k_0)\in\psi_2$. Then we have $\er(\gamma_\delta)=\er(\gamma_0)$ and 
$\gamma_\delta=:\gamma_1$ does dot depend on $\delta$. Then we have $(\er(\delta);s(\gamma_1),\gamma_0)\in m'(\phi\times\phi)=\phi m$, so 
$(s(\gamma_1),\gamma_0)\in D(m)$ and the element $s(\gamma_1) \gamma_0$ is in relation with each element $\er(\delta)$ for 
$\delta\in \phi(\gamma_0)$, but these are all 
elements in $\Delta$ that $s(\gamma_1) \gamma_0$ can be in relation with, so this element is in $\ker(\phi)$.

\noindent
$\Rightarrow$  Now suppose $\ker(\phi)\neq \Gamma^{(0)}$. One reason for that can be inequality
$D(\phi)\neq \G$. Since $D(\phi)$ is a union of orbits there exists an orbit $O$ and
$O\cap D(\phi)=\emptyset$. Let $\G_1:=\G^{(0)}$ (set-groupoid) and  define  mappings $f_1,f_2:\G^{(0)}\rightarrow\G^{(0)}$ by
$f_1:=id$ and  $f_2$ by $f_2(e)=e$ for $e\in \G^{(0)}\setminus O$ and $f_2(e)=e_0$ for
$e\in O$ and some $e_0\in \G^{(0)}\setminus O$. For $h_1:=f_1^T$, $h_2:=f_2^T$ we have 
$\phi h_1=\phi h_2$ but $h_1\neq h_2$ so $\phi$ is not a monomorphism. 

Let $D(\phi)=\G^{(0)}$ but $\ker(\phi)\neq \G^{(0)}$, so there exists $\gamma_0\not\in \G^{(0)}$ and $\gamma_0\in \ker(\phi)$.
Let $e_0=\el(\gamma_0)=\er(\gamma_0)$, $G_0$ be the isotropy group of $e_0$ and $H_0:=G_0\cap \ker(\phi)$ -- this is a  subgroup of $G_0$.
Consider two relations $\psi_1,\psi_2:H_0\rel\Gamma$ defined by $\psi_1:=\G^{(0)} \times H_0$ and 
$\psi_2:=(\G^{(0)}\setminus\{e_0\})\times H_0\cup diag(H_0\times H_0)$
Then $\psi_1,\psi_2$ are morphisms $H_0\rel G$, $\psi_1\neq \psi_2$ but
$\phi \psi_1=\phi \psi_2=\Delta^{(0)}\times H_0$.\\
\dowl

\begin{col} Every morphism $h:\Gamma\rel\Delta$ of groupoids defines a homomorphism $\tilde{h}$ of groups of bisections (lemma \ref{mor-and-bis}). 
The previous proposition implies that $h$ is mono iff $\tilde{h}$ is mono (i.e. is injective).
\end{col}

\subsection{ Examples of monomorphisms}

Any morphism $h: X^2 \rel \Delta$ is a monomorphism 
(because $X^2$ is transitive and has trivial isotropy).

If $\Gamma_1\subset \Gamma$ is a wide subgroupoid, the inclusion $i:\Gamma_1\rightarrow\Gamma$ is a monomorphism 
(i.e. wide subgroupoids are subobjects).

For a groupoid $\Gamma$ the left regular representation $l:\Gamma\rel \Gamma^2$ is a monomorphism.

For a groupoid $\Gamma\rightrightarrows E$ the mapping $\G\ni \gamma\mapsto(e_L(\gamma),e_R(\gamma))\in E^2$ is monomorphism iff 
$\G$ has trivial isotropy (the kernel of this morphism is just the isotropy group bundle).

For a bisection $B\subset\G$ the mapping $Ad_B:\Gamma\rightarrow\Gamma$ is a monomorphism.

If a group $G$ acts on $X$ and $h:G\rel X^2$ is the associated morphism, then $h$ is monomorphism iff the action is effective 
(i.e. only neutral element of $G$ acts as identity mapping).

\subsection{Epimorphisms}
In the category of groups epimorphisms are surjective group homomorphisms. In groupoids we have:
\begin{prop}
Let $h:\Gamma\rel \Delta$ be an epimorphism. Then $h$ is surjective (i.e. $Im(h)=\Delta$).
\end{prop}
{\em Proof:} The proof is an adaptation of the proof for groups; 
since an image of a  morphism is a wide subgroupoid it is enough to prove: 
{\em Let  $\Gamma_1$ be a wide subgroupoid of $\Gamma$ and $\Gamma_1\neq \Gamma$. There exists a  groupoid $\Gamma_2$ and morphisms 
$h_1,h_2:\Gamma\rel\Gamma_2$ such that $h_1\neq  h_2$ but $h_1$ and $h_2$ coincide on $\Gamma_1$}

Let $\gamma_0\in\G\setminus\G_1$, $e_0:=\el(\gamma_0)$ and  $H:=\er^{-1}(e_0)\cap \G_1$. Since $\G_1$ is a  subgroupoid, sets 
$H$ and $H\gamma_0$ are disjoint and the right multiplication by $\gamma_0$ is a bijection $H\rightarrow H\gamma_0$.
Define an involution $\sigma:\Gamma\rightarrow \Gamma$ by:
$$\sigma(\gamma):=\left\{\begin{array}{lcr} \gamma \gamma_0 &for & \gamma\in H\\\gamma s(\gamma_0) & for &\gamma\in H\gamma_0\\
\gamma& & otherwise\end{array}\right.$$
Note that $(\gamma_1,\gamma_2)\in D(m)\iff (\gamma_1,\sigma(\gamma_2))\in D(m)$.
 
\noindent Let $\tilde{\sigma}$ be the bisection of $\G^2$ defined by $\sigma$ (i.e. $\tilde{\sigma}:=\{(\sigma(\gamma),\gamma): \gamma\in \Gamma\}$).
We claim that  morphisms $l$ and $Ad_{\tilde{\sigma}}l: \G\rel\G^2$  coincide on $\G_1$.
By the definition of $l$ (formula (\ref{left})) and definition of $\tilde{\sigma}$ we have:
$$(\gamma_1,\gamma_2;\gamma_3)\in l\iff (\gamma_1;\gamma_3,\gamma_2)\in m\,\,{\rm  and} \,\, 
(\gamma_1,\gamma_2;\gamma_3)\in Ad_{\tilde{\sigma}}l\iff (\sigma(\gamma_1);\gamma_3,\sigma(\gamma_2))\in m.$$
For $\gamma_3\in\Gamma_1$ we have:\\
if $\gamma_2\in H$  and $\gamma_1=\gamma_3 \gamma_2$ then $\gamma_1\in H$ so 
$\sigma(\gamma_1)=\gamma_1\gamma_0=\gamma_3\gamma_2\gamma_0=\gamma_3\sigma(\gamma_2)$;\\
if $\gamma_2\in H\gamma_0$  and $\gamma_1=\gamma_3 \gamma_2$ then $\gamma_1\in H\gamma_0$ so 
$\sigma(\gamma_1)=\gamma_1 s(\gamma_0)=\gamma_3\gamma_2s(\gamma_0)=\gamma_3\sigma(\gamma_2)$;\\
if $\gamma_2\not\in (H\cup H\gamma_0)$ then $\gamma_1\not\in (H\cup H\gamma_0)$ so 
$\sigma(\gamma_1)=\gamma_1=\gamma_3\gamma_2=\gamma_3\sigma(\gamma_2)$.\\
In this way for $\gamma_3\in\Gamma_1: \,(\gamma_1,\gamma_2;\gamma_3)\in l \Rightarrow (\gamma_1,\gamma_2;\gamma_3)\in Ad_{\tilde{\sigma}}l$.
The opposite inclusion follows if we substitute $\sigma(\gamma_1)$ and $\sigma(\gamma_2)$ instead of $\gamma_1$ and $\gamma_2$.

It remains  to check whether $l\neq Ad_{\tilde{\sigma}}l$. 
Certainly we have $(e_0;s(\gamma_0),\gamma_0)\in l$ and this element is contained in 
$Ad_{\tilde{\sigma}}l$ iff $\gamma_0=\gamma_0 \sigma(s(\gamma_0))$ i.e.
$\er(\gamma_0)=\sigma(s(\gamma_0))$. Since $s(\gamma_0)$ is neither in $H$ nor  in $E$, this equality implies that 
$s(\gamma_0)=\gamma'\gamma_0$ for some $\gamma'\in H$ but then $\sigma(s(\gamma_0))=\gamma'=\er(\gamma_0)$ i.e. 
$\el(\gamma_0)=\er(\gamma_0)$ and $\gamma_0=s(\gamma_0)$.
In this way the proposition is proven if there exists $\gamma_0\in\G\setminus\G_1$ such that $s(\gamma_0)\neq \gamma_0$.

Let us now assume that every element $\gamma\in \Gamma\setminus \G_1$ satisfies $s(\gamma)=\gamma$. Choose some 
$\gamma_0\in\Gamma\setminus \G_1$, we claim that orbit of $e_0:=\er(\gamma_0)=\el(\gamma_0)$ is equal to $\{e_0\}$. 
Indeed if $\er(\gamma)=e_0$ and $\el(\gamma)=e_1\neq e_0$ then $\gamma\in\Gamma_1$, but then $\gamma\gamma_0\in\G\setminus\G_1$ 
and this is impossible. So there exists one point orbit $\{e_0\}$ in $E$ i.e $\Gamma$ is a union of a groupoid $\tilde{\Gamma}$ and 
$G$ - the isotropy group of $e_0$.  Then the intersection $H:=\Gamma_1\cap G$  ia a normal subgroup of $G$. 
Let $K:=G/H$ and $\pi: G\rightarrow K$ be the canonical homomorphism. Let $e_K$ be the neutral element in $K$ and 
define two morphisms $\Gamma\rel K$:  $h_1:=\{(e_K,g); g\in G\}$  and $h_2:=\{(\pi(g),g): g\in G\}$. 
It is clear that $h_1\neq h_2$ but they coincide on $\Gamma_1$.\\
\dow

\noindent Being surjective is  not sufficient for being an epimorphism:
\begin{ex1}
Let a group $G$ act  transitively on a set $X$ and assume that $\phi(g x)=g\phi(x)$ 
for some bijection $\phi:X\rightarrow X$. Let $h:G\rel X^2$ be the associated morphism: $Gr(h):=\{(gx,x;g):g\in G\,,\,x\in X\}$ 
and $\tilde{\phi}$ the bisection of $X^2$ defined by $\phi$. Then  straightforward computation shows that $h=Ad_{\tilde{\phi}} h$.
\end{ex1}
If  $h:\G\rel \G_1$ is a morphism which is a surjective mapping $\G\supset D(h)\rightarrow \G_1$ then $h$ is 
an epimorphisms (such mappings are epimorphisms in the category of sets and relations), but being a mapping is not necessary:
\begin{prop}
Let a group $G$ act  transitively on a set $X$ and assume that for every $x$ there exists $g$ such that 
$x$ is {\em the unique} fixed point for $g$.
Then the  morphism associated to this action is an epimorphism $G\rel X^2$.
\end{prop}
{\em Proof: }
Assume that  $h_1,h_2:X^2\rel\Gamma$ are morphisms such that  $h_1 h=h_2 h$ and let  $\rho_1,\rho_2$ be their base maps.
We will show  that $\rho_1=\rho_2$.
Indeed, if $x=\rho_1(e)$ then, by the assumption, there is $g$ such that $(e,g)\in h_1 h$ and $x$ is the unique fixed point of $g$. 
Since $h_1 h=h_2 h$  there exist $(x_1,x_2)$ such that $(x_1, x_2; g)\in h$ and $(e; x_1, x_2)\in h_2$. 
But then $x_1=x_2$ and $x_1=g x_1$ so $x_1=x$ and $x=\rho_2(e)$.
The opposite inclusion can be shown in the same way.
The following lemma states this is sufficient:

\begin{lem} Let $h:\Gamma\rel\Delta$ be a surjective morphism, $h_1,h_2:\Delta\rel \Lambda$ morphisms with base maps 
$\rho_1, \rho_2:\Lambda^{(0)}\rightarrow \Delta^{(0)}$. If $h_1 h=h_2 h$ and $\rho_1=\rho_2$ then $h_1=h_2$.
\end{lem}

{\em Proof:} Let $(\lambda,\delta)\in h_1$; exists $\gamma\in\Gamma$ such that $(\delta,\gamma)\in h$ so $(\lambda,\gamma)\in h_1 h=h_2 h$, so  
exists $\delta_1$ such that: $(\delta_1,\gamma)\in h$ and $(\lambda,\delta_1)\in h_2$. But since $\rho_1=\rho_2$ it follows that 
$\er(\delta)=\er(\delta_1)$, so $\delta=\delta_1$. In this way $Gr(h_1)\subset Gr(h_2)$ and the opposite inclusion can be proved in the same way.

\dow

\subsection{Examples of epimorphisms}

If $\G_1\subset\G$ is a union of transitive components and $i:\G_1\rightarrow \G$ the inclusion map, then $i^T:\G\rel\G_1$ is an epimorphism.

Let $\Gamma\rra E$ be a groupoid and $R\subset E\times E$ its orbit relation. The map 
$\pi: \Gamma\ni\gamma\mapsto (\el(\gamma),\er(\gamma))\in R$ is an epimorphism.

Let $\Gamma\rra E$ be a groupoid, $\G'$ its isotropy group bundle and $R\subset E\times E$ its orbit relation. Let 
$j:\G'\rel \G$ be inclusion and $\pi: \Gamma\ni\gamma\mapsto (\el(\gamma),\er(\gamma))\in E\times E$. Then $\pi:\G\rel R$ is an 
epimorphism, $j:\G'\rel\G$ is a monomorphism and $ker(\pi)=Im(j)$. {\em In other words every groupoid is an extension of equivalence 
relation by a group bundle.}

If $G=Bij(X)$ is the group of bijection of a set $X$ then the morphism $G\rel X^2$ coming from canonical (left) action of $G$ on $X$
is mono- and epimorphism but not an isomorphism.

If $\Gamma\neq \Gamma^{(0)}$ then the left regular representation $l:\Gamma\rel \G^2$ is not an epimorphisms. In fact, in this situation there
 exists a bisection $B\neq \Gamma^{(0)}$. Let $R$ be the right multiplication by $B$ and $\tilde{R}$ the corresponding bisection of $\G^2$,
then $l=Ad_{\tilde{R}} l$.

If $\Gamma=\Gamma_1\cup\Gamma_2$ is a disjoint union then ``the canonical projections'' $i_1^T(i_2^T):\Gamma\rel \Gamma_1 (\Gamma_2)$ 
are epimorphisms. With this ``projections'' we have:
\begin{prop} Disjoint union is a product in the category of groupoids.
\end{prop}\dow

\section{Groupoid actions and morphisms}
A groupoid $\Gamma\rightrightarrows E$ can act on a set $X$ equipped with a mapping to $E$. 
The standard definition is:
\begin{defi}\label{action-standard}
Let $\G\rightrightarrows E$ be a groupoid, $X$ a set and $\rho:X\rightarrow E$ a mapping. 
Define the set  $\G \vphantom{x}_{e_R}\hspace{-1ex}\times_{\rho} X:=\{(\gamma,x)\in\G\times X: e_R(\gamma)=\rho(x)\}$.
An action of $\G$ on $X$ is a mapping: $\G \vphantom{x}_{e_R}\hspace{-1ex}\times_{\rho} X\ni(\gamma,x)\mapsto \gamma x\in X$ that satisfies:
$\rho(x) x=x$ and $\gamma_1(\gamma_2 x)=(\gamma_1 \gamma_2) x$ i.e. if one side is defined, the other is also and then they are equal.
\end{defi}
\noindent If $\Gamma$ is a group this definition gives just a group action. Basic examples for general groupoids  are :

\begin{ex}
Action by multiplication (from the left) : $X:=\Gamma\,,\,\rho:=\el\,,\,\G \vphantom{x}_{e_R}\hspace{-1ex}\times_{\rho} X=\Gamma^{(2)}$ 
and the action is $(\gamma,\gamma')\mapsto \gamma \gamma'$.

Action on the set of units:  $X:=E\,,\,\rho=id\,,\,\G \vphantom{x}_{e_R}\hspace{-1ex}\times_{\rho} X=\{(\gamma,\er(\gamma))\,,\,\gamma\in\Gamma\}$
and the action is: $(\gamma,\er(\gamma))\mapsto \el(\gamma)$

 Action on the isotropy group bundle: $X=\G'$ -- the isotropy group bundle of $\G$; $\rho=\el\,,\,
\G \vphantom{x}_{e_R}\hspace{-1ex}\times_{\rho} X=\Gamma^{(2)}\cap(\G\times\G')$ and the action is $(\gamma,\gamma')\mapsto \gamma \gamma' s(\gamma)$.
\end{ex}

If we use relations the definition of a groupoid action can be presented in a more group-like style:
\begin{defi}\label{action}
Let $\Gamma\rightrightarrows E$ be a groupoid and $X$ a set. An action of $\Gamma$ on $X$ is a
 relation $\Phi:\Gamma\times X\rel X$ that satisfies: $$\Phi(m\times id)=\Phi(id\times \Phi)\,,\,\,\Phi(e\times id)=id.$$
\end{defi}
\noindent Next proposition states the equivalence of both definitions.
\begin{prop}\label{action-equiv}
Let $\Phi: \Gamma\times X\rel X$ be an action in a sense of def. \ref{action}. Then
\begin{enumerate}
\item For every $x\in X$ there exists unique $e\in E$ such that $(x;e;x)\in \Phi$, i.e. $\Phi$ defines a mapping
$\rho:X\rightarrow E$;
\item $D(\Phi)=\G \vphantom{x}_{e_R}\hspace{-1ex}\times_{\rho} X$; 
\item $(y;\gamma, x)\in\Phi\Rightarrow \rho(y)=\el(\gamma)$;
\item $(y;\gamma,x)\in \Phi\iff (x;s(\gamma),y)\in \Phi$;
\item $\Phi$ is a mapping $D(\Phi)\rightarrow X$;  this mapping is an action of $\Gamma$  on $X$ in the sense of def. \ref{action-standard};
\item If $\Gamma$ acts on $X$ in a sense of def. \ref{action-standard}, the relation $\Phi:=\{(\gamma x; \gamma, x): \er(\gamma)=\rho(x)\}$ 
is an action in the sense of def. \ref{action}.
\end{enumerate}
\end{prop}
{\em Proof:} Straightforward computations give:
$$\Phi(m\times id)=\{(y;\gamma_1,\gamma_2,x): (\gamma_1,\gamma_2)\in D(m)\,,\,(y;\gamma_1\gamma_2,x)\in \Phi\}$$
$$\Phi(id\times \Phi)=\{(y,\gamma_1,\gamma_2,x): \exists z: (y;\gamma_1,z)\in\Phi \,,\,(z;\gamma_2,x)\in \Phi\}$$
1) The condition $\Phi(e\times id)=id$ means that for any $x\in X$ there exists $e\in E$ such that $(x;e,x)\in \Phi$. 
Suppose that $(x;e_1,x)\in\Phi$ and $(x;e_2,x)\in\Phi$ i.e $(x;e_1,e_2,x)\in \Phi(id\times \Phi)=\Phi(m\times id)$. 
Therefore  $(e_1,e_2)\in D(m)$ so  $e_1=e_2$ and  the first statement is proven.\\
2) For any $x\in X$, $(x;\rho(x), x)\in \Phi$; 
let $e_R(\gamma)=\rho(x)$ so $(\rho(x);s(\gamma),\gamma)\in m$ and $(x;s(\gamma),\gamma,x)\in \Phi(m\times id)=\Phi(id\times \Phi)$. 
Therefore exists $z$ such that: $(x;s(\gamma),z)\in \Phi$ and $(z;\gamma,x)\in \Phi$, i.e $(\gamma,x)\in D(\Phi)$. \\
Conversely, if $(\gamma,x)\in D(\Phi)$ i.e. $(z;\gamma,x)\in\Phi$ for some $z\in X$, then (since $(\gamma;\gamma,e_R(\gamma))\in m$) 
$(z;\gamma,e_R(\gamma),x)\in \Phi(m\times id)=\Phi(id\times\Phi)$, so there exists $y$ such that $(z;\gamma,y)\in\Phi$ and
$(y;e_R(\gamma),x)\in \Phi$. It follows that $y=x$ and $e_R(\gamma)=\rho(x)$.\\
3) Let $(y;\gamma,x)\in\Phi$. Since $(y;\rho(y), y)\in \Phi$ we have $(y;\rho(y),\gamma,x)\in \Phi(id\times\Phi)=\Phi(m\times id)$, i.e. 
exists $\gamma_1$ such that: $(y;\gamma_1,x)\in \Phi$ and $(\gamma_1;\rho(y), \gamma)\in m$. The last condition implies $\rho(y)=\el(\gamma)$.\\
4) Since $\gamma=s(s(\gamma))$ it is enough to prove implication: $(y;\gamma,x)\in \Phi\Rightarrow (x;s(\gamma);y) \in \Phi$.
If $(y;\gamma,x)\in\Phi$ then (by the previous statement)  exists $z\in X$ such that 
$(z;s(\gamma),y)\in\Phi$ so $(z;s(\gamma),\gamma,x)\in\Phi(id \times \Phi)=\Phi(m\times id)$. 
It follows that $(z;\er(\gamma),x)\in \Phi$ i.e. $z=x$.\\  
5) Let $(y_1;\gamma,x)\in\Phi$ and $(y_2;\gamma, x)\in \Phi$. Then (by the previous statement) $(x;s(\gamma),y_2)\in\Phi$ and 
$(y_1;\gamma,s(\gamma),y_2)\in \Phi(id\times\Phi)=\Phi(m\times id)$; so $(y_1;\el(\gamma),y_2)\in \Phi$ and $y_1=y_2$; i.e. 
$\Phi:\G \vphantom{x}_{e_R}\hspace{-1ex}\times_{\rho} X\rightarrow X$ is a mapping. 
It is straightforward that this mapping satisfies conditions in  def. \ref{action-standard}.\\
6) The last statement is also straightforward.

\dow

\noindent
A morphism from a group $G$ to a pair groupoid $X^2$ defines an action of $G$ on $X$ (see examples in section 3);
the similar fact  is  true for any groupoid:
\begin{lem}
Let $X$ be a set, $\G$ a groupoid and $h:\Gamma\rel X^2$ a morphism. The relation $\Phi: \Gamma\times X\rel X$ defined by
$$(x_1;\gamma, x_2)\in\Phi\iff (x_1,x_2;\gamma)\in h$$ is an action of $\Gamma$ on $X$. Conversely, if  $\Phi:\Gamma\times X\rel X$ is an action 
then $h:=\{(\Phi(\gamma,x),x;\gamma): (\gamma,x)\in D(\Phi)\}$ is a morphism $\Gamma\rel X^2$.\\
\dowl
\end{lem}

Any morphism $h: \Gamma\rel\Delta$, by a composition with $l_\Delta: \Delta\rel \Delta^2$ (the left regular representation 
of $\Delta$),  defines a morphism 
$l_\Delta h:\Gamma\rel\Delta^2$ i.e an action of $\Gamma$ on $\Delta$. This action commutes with multiplication (from the right) in $\Delta$.
Moreover it  turns out that morphisms of groupoids are precisely those actions that commutes with multiplication (from the right):
\begin{prop}
Let $\Gamma,\Delta$ be groupoids and $\Phi:\Gamma\times \Delta\rel \Delta$ an action. If 
$$\Phi(id\times m_\Delta)=m_\Delta(\Phi\times id)$$ 
then $h:\Gamma\rel \Delta$ defined by $h:=\{(\Phi(\gamma,\delta) s(\delta) ;\gamma): (\gamma,\delta)\in D(\Phi)\}$ is a morphism.
\end{prop}
{\em Proof:} The assumed equality of compositions is equivalent to:
$$\left[(\gamma,\delta_2)\in D(\Phi)\,,\,(\delta_1;\Phi(\gamma,\delta_2),\delta_3)\in m\right]\iff
\left[(\delta_2,\delta_3)\in D(m)\,,\,(\delta_1;\gamma,\delta_2\delta_3)\in\Phi\right]$$
So for any $(\gamma,\delta)\in D(\Phi)$ we have $(\gamma,\el(\delta)=\delta s(\delta))\in D(\Phi)$ 
and $h=\{(\Phi(\gamma,\el(\delta)),\gamma): (\gamma,\delta)\in D(\Phi)\}$. It is rather straightforward that $h$ is a morphism. \cite{BunPs}\\
\dow

\subsection{Morphisms as functors between action categories}
{\em This  point of view on morphisms was communicated to the author by Ralf Meyer.} 
\begin{defi}
Let $\G$ be a groupoid; a $\G$-set is a pair $(X,\Phi)$, where $X$ is a set and $\Phi$ an action of $\G$ on $X$. 
Let  $(X,\Phi)$ and $(Y,\Psi)$ be $\G$-sets. A map $f:X\rightarrow Y$ is equivariant iff $f\Phi=\Psi(id\times f)$.
\end{defi}
\begin{re} In the classical book of Renault \cite{Ren}  ``$\Gamma$-sets'' mean bisections; but since ``bisections''  
seem to be widely accepted we can use $\Gamma$-sets as in the theory of group actions.
\end{re}
$\G$-sets with  equvariant maps as morphisms form a category. 
If we think of actions as of morphisms to pair groupoids, an  equivariant map $f:X\rightarrow Y$ is characterized by
$(f\times id)h_1=(id\times f^T)h_2$ for $h_1:\Gamma\rel X^2$ and $h_2:\Gamma\rel Y^2$.
A morphism $h:\G\rel \Delta$ defines a functor $H_h$ from $\Delta$-sets to $\G$-sets by
composition. This functor doesn't change  sets and equivariant maps, in other words, if $For_\Gamma, For_\Delta$ are forgetful functors 
to the category of sets (i.e $For_\Gamma(X,\Phi)=X$ and $For_\Gamma(f)=f$, where
$f$ is an equivariant map between $\Gamma$-sets $X$ and $Y$) it satisfies  $For_\Delta H_h=For_\Gamma$. 
Conversely any such functor defines  a morphism of groupoids:
\begin{prop} Let $H$ be a functor from $\Gamma$-sets to $\Delta$-sets satisfying
$For_\Delta H=For_\Gamma$. There exists unique morphism $h:\Delta\rel\Gamma$, such that $H$ is the composition with $h$.
\end{prop}
{\em Proof:} Let $E$ denote the set of units of $\Gamma$. For $e\in E$ let $\Gamma_e:=\er^{-1}(e)$ and consider $\Gamma$-set $(\Gamma_e, m_e)$, where 
 the action  $m_e: \Gamma\times\Gamma_e\rel\Gamma_e$ is the  multiplication. Applying the functor $H$ we get $\Delta$-set 
$(\Gamma_e,H(m_e))$, where $H(m_e):\Delta\times \Gamma_e\rel \Gamma_e$ is an action of $\Delta$. 
Define a relation $\Phi:\Delta\times\Gamma\rel\Gamma$ by 
$$(\gamma_1;\delta,\gamma)\in \Phi\iff (\gamma_1;\delta,\gamma)\in H(m_e)\,\,{\rm for}\,\, e:=e_R(\gamma)$$
It is rather straightforward to verify that $\Phi$ is an action; let us check it commutes with right multiplication in $\Gamma$ i.e. 
defines a  morphism from $\Delta$ to $\Gamma$. 
$$(\gamma_1;\delta_1,\gamma_2,\gamma_3)\in \Phi(id\times m_\Gamma)\iff 
(\gamma_2,\gamma_3)\in D(m_\Gamma)\,,\, (\gamma_1;\delta_1,\gamma_2\gamma_3)\in\Phi$$
what is equivalent to  $\er(\gamma_2)=\el(\gamma_3)$ and $(\gamma_1;\delta_1,\gamma_2 \gamma_3)\in H(m_e)$ for $e:=\er(\gamma_3)$.
On the other hand:
$$(\gamma_1;\delta_1,\gamma_2,\gamma_3)\in m_\Gamma(\Phi\times id)\iff (\delta_1,\gamma_2)\in D(\Phi)\,,\,
\gamma_1=\Phi(\delta_1,\gamma_2)\gamma_3$$
what is equivalent to $\er(\gamma_1)=\er(\gamma_3)$ and $(\gamma_1 s(\gamma_3);\delta, \gamma_2)\in H(m_{\tilde{e}})$ for $\tilde{e}:=\er(\gamma_2)$.

So let $(\gamma_1;\delta_1,\gamma_2,\gamma_3)\in \Phi(id\times m_\Gamma)$; $e:=\er(\gamma_3)$, $\tilde{e}:=\el(\gamma_3)=\er(\gamma_2)$. 
Let $R$ be the  right multiplication by $s(\gamma_3)$ i.e. $R:\er^{-1}(e)\ni\gamma\mapsto \gamma s(\gamma_3)\in \er^{-1}(\tilde{e})$. $R$ is a morphism
from $(\Gamma_e, m_e)$ to $(\Gamma_{\tilde{e}}, m_{\tilde{e}})$: $R m_e=m_{\tilde{e}}(id \times R)$. Applying the functor $H$ we get $H(R)=R$ 
(by assumption on $H$) and 
$RH(m_e)=H(m_{\tilde{e}})(id \times R)$. If $(\gamma_1;\delta_1,\gamma_2 \gamma_3)\in H(m_e)$ then $(\gamma_1 s(\gamma_3);\delta_1,\gamma_2 \gamma_3)\in R H(m_e)=H(m_{\tilde{e}})(id \times R)$, i.e. $(\gamma_1 s(\gamma_3);\delta_1,\gamma_2)\in H(m_{\tilde{e}})$. 
The opposite inclusion can be shown in the same way.\\
\dow

\subsection{Action groupoids and functors}
Let $\Phi:\Gamma\times X\rel X$ be an action of a groupoid $\Gamma\rightrightarrows E$ on  
$X$ with a base map $\rho:X\rightarrow E$. 
It follows from prop \ref{action-equiv} that the following definitions make sense:
$$E_\Phi:=\{(\rho(x),x):x\in X\}\,,\,\,s_\Phi: D(\Phi)\ni(\gamma,x)\mapsto (s(\gamma),\Phi(\gamma,x))\in D(\Phi)$$
$$m_\Phi:D(\Phi)\times D(\Phi)\rel D(\Phi)\,,\,
Gr(m_\Phi):=\{(\gamma_1 \gamma_2,x;\gamma_1,\Phi(\gamma_2,x),\gamma_2,x):(\gamma_1,\gamma_2)\in D(m)\,,\,(\gamma_2,x)\in D(\Phi)\}$$
\begin{lem} $(D(\Phi),m_\Phi,s_\Phi,E_\Phi)$ is a groupoid; it is called  the action groupoid for an action $\Phi$ and 
will be denoted by $\Gamma\times_\Phi X$.\\\dowl\end{lem}

Let $\Gamma\rightrightarrows E$ and $\Delta \rightrightarrows F$ be groupoids and $h:\Gamma\rel\Delta$ a morphism. Composition of $h$ 
with the mapping $\Delta\ni\delta\mapsto(\el(\delta),\er(\delta))\in F^2$ gives a morphism $\Gamma\rel F^2$, i.e. the action 
$\phi_h: \Gamma\times F\rel F$.  Its domain is  $D(\phi_h):=\{(\gamma,f):\er(\gamma)=\rho_h(f)\}$ and the action is
$(\gamma,f)\mapsto \el(h^R_f(\gamma))$. The relation $h$ defines a mapping $D(\phi_h)\ni(\gamma,f)\mapsto h^R_f(\gamma)\in\Delta$.
This mapping is a functor from the action groupoid $\Gamma\times_{\phi_h} F$ to $\Delta$.  

Conversely: an action $\phi$ of $\Gamma$ on $F$ and a functor $K:\Gamma\times_\phi F\rightarrow \Delta$ satysfying
$K(e,f)=f$ for $(e,f)\in D(\phi)\cap (E\times F)$ defines a morphism $h:\Gamma\rel \Delta$ by 
$Gr(h):=\{(K(\gamma,f),\gamma): (\gamma,f)\in \Gamma\times_\phi F\}$ (this was observed in  \cite{Bun}).

\subsection{``Homogenous spaces'' and quotient groupoids}
Let $\G\rightrightarrows E$ be a groupoid and $G\subset \G$ a wide subgroupoid (i.e. $E\subset G$); consider a relation on $\G$:
$\gamma_1\sim\gamma_2 \iff s(\gamma_1) \gamma_2\in G$.
It is straightforward to check that this is an equivalence relation. Let $\G/G$ denotes the set of equivalence classes, 
$\pi:\G\rightarrow \G/G$ be the canonical projection and define a relation $\tilde{m}_G:\G\times \G/G\rel \G/G$ by:
$$Gr(\tilde{m}_G):=(\pi\times id\times \pi)(Gr(m)).$$
\begin{prop}
$\tilde{m}_G$ is an action of $\G$ on $\G/G$.\\
\vspace{-1ex}\dow
\end{prop} 
\begin{ex}\label{ex-tran}
For $G=\Gamma$ we get $\Gamma/\Gamma= E$,  $\pi=\el$ and the action is just the  action of $\Gamma$ on its set of units as in the 
beginning of this section: $(\gamma,\er(\gamma))\mapsto \el(\gamma)$;

For $G=\Gamma'$($=$the isotropy group bundle of $\Gamma$), $\Gamma/\Gamma'$ can be identified with the orbit equivalence relation 
$R\subset E\times E$ by the map $\Gamma/\Gamma' \ni[\gamma]\mapsto (\el(\gamma),\er(\gamma))\in R$ and the action $\tilde{m}_G:\G\times R\rel R$ is 
$\{(\el(\gamma),e;\gamma,\er(\gamma),e): e\in E, \gamma\in\Gamma\}$.

Let $G=R\subset X\times X$ be an equivalence relation, it is a wide subgroupoid of $X^2$; let $Y:=X/\sim$ be the set of classes. Then
$X^2/R$ may be identified with $X\times Y$ and the action is: 
$\{(x_1,y;x_1,x_2,x_2,y): x_1,x_2\in X\,,\,y\in Y\}\subset (X\times Y)\times X^2\times (X\times Y)$

\end{ex}
\noindent If $G\subset \G$ satisfies some additional conditions, it is possible to define a groupoid structure on $\G/G$:
\begin{prop}
Let a wide subgroupoid $G\subset\G$ be  contained in the  isotropy group bundle and (for any $e\in E$) 
$G\cap \G_e$
be a normal subgroup of $\G_e$.
Let $\pi:\G\rightarrow \G/G$ be the canonical projection. 
Define $m_G:\G/G\times \G/G\rel \G/G$ by $Gr(m_G):=(id\times \pi\times id)(Gr(\tilde{m}_G))$ and  $E_G:=\pi(E)$.
\begin{itemize}
\item $\pi(s(\gamma))$ depends only on $\pi(\gamma)$ i.e. formula $s_G(\pi(\gamma)):=\pi(s(\gamma))$ defines a mapping $s_G:\G/G\rightarrow\G/G$;
\item $(\G/G, m_G,s_G,E_G)$ is a groupoid;
\item $\pi:\G\rightarrow \G/G$ is an epimorphism.
\end{itemize}
\vspace{-1ex}\dow
\end{prop}
If a domain of morphism is whole groupoid, its kernel is a wide subgroupoid that satisfies assumptions of previous proposition. 
As for groups we can then factorize it:
\begin{lem}
Let $h:\Gamma\rel\Gamma_1$ be a morphism with   $D(h)=\Gamma$. Let $\tilde{\Gamma}:=\Gamma/\ker(h)$ and $\pi:\Gamma\rightarrow \tilde{\Gamma}$ 
be the canonical projection. The relation $\tilde{h}:\tilde{\Gamma}\rel\Gamma_1$ defined by $Gr(\tilde{h}):=(id\times \pi) Gr(h)$ is a monomorphism and
$h=\tilde{h}\cdot \pi$.
\end{lem}
{\em Proof:} Let $\gamma_0\in\ker(h)$ and $(\gamma,\gamma_0)\in D(m)$. If $(\delta,\gamma)\in h$ then $(\er(\delta),\gamma_0)\in h$ and
$(\delta;\gamma,\gamma_0)\in m'(h\times h)=h m$ so $(\delta;\gamma \gamma_0)\in h$. In other words $\tilde{h}$ is defined by:
$$(\delta,[\gamma])\in\tilde{h}\iff (\delta,\gamma)\in h.$$
It is routine to check that $\tilde{h}$ is a morphism and monomorphism.
\dowl
 
Domain $D(h)$ of a morphism $h:\Gamma\rel\Gamma_1$ is a union of transitive components, so  $i^T:\G\rel D(h) $ is epimorphism, 
where $i:D(h)\rightarrow \G$ is the inclusion map. We can restrict $h$ to $D(h)$ i.e. consider a morphism $h_1:D(h)\rel \Gamma_1$, then 
$h=h_1\cdot i^T$ is a composition of morphisms. Using the previous lemma one can prove
the following factorisation of morphisms:
\begin{prop}\label{factor}
Let $h:\Gamma\rel \Delta$ be a morphism. There exist a groupoid $\tilde{\Gamma}$, epimorphism $h_1:\Gamma\rel \tilde{\Gamma} $ 
and monomorphism $h_2:\tilde{\Gamma}\rel \Delta$ such that $h=h_2 h_1$.
\end{prop}
\vspace{-1ex}\dow
\subsection{Transitivity}

For a group $G$ its  transitive actions are actions on homogenous spaces, i.e. if $\Phi:G\times X\rightarrow X$ is a transitive action, 
there exists a subgroup $H\subset G$ and a bijection  $\psi: X\rightarrow G/H$ such that $\tilde{m}_H(id\times \psi)=\psi\Phi $, 
where $\tilde{m}_H(g,g'H):=g g' H$.

Let $\Phi:\G\times X\rel X$ be an action of groupoid $\G$ and $\rho_\Phi$ its base map. 
If $\rho_\Phi(x_1)$ and $\rho_\Phi(x_2)$ are in different orbits, 
there is no $\gamma\in\G$ such that $(x_1;\gamma,x_2)\in\Phi$. So if $\G$ is not transitive, even the action $m:\G\times \G\rel \G$ is  
{\em not transitive} in this sense, moreover  the third example in Ex. \ref{ex-tran} shows that even for transitive groupoids an 
action on ``homogenous space''  need not be transitive. The following proposition tells how   to identify homogenous spaces for transitive groupoids.
Note that if $\G$ is a group i.e. $E=\{e\}$ -- the neutral element of $\G$, the condition $\Phi(\Gamma\times p(E))=X$ used below means 
transitivity of an action.
\begin{prop}
Let $\G\rra E$ be a transitive groupoid and $\Phi: \G\times X\rel X$ an action with base map $\rho:X\rightarrow E$.
Assume  there exists $p:E\rightarrow X$ which is a section of $\rho$ such that $\Phi(\Gamma\times p(E))=X$. Then 
there exists (wide) subgroupoid $G\subset \G$ and a 
bijection $\psi: X\rightarrow \G/G$ such that $\tilde{m}_G(id\times \psi)=\psi\Phi $
\end{prop}
{\em Proof:}   For $(e_1,e_2)\in E\times E$ define 
$G(e_1,e_2):=\{\gamma\in\Gamma: (p(e_1);\gamma,p(e_2))\in\Phi\}$;
basic properties of these sets are:  
$$G(e_1,e_2)\subset \el^{-1}(e_1)\cap \er^{-1}(e_2)\,;\,\,\forall e\in E: e\in G(e,e)\,;\,\,\gamma\in G(e_1,e_2)\Rightarrow s(\gamma)\in G(e_2,e_1);$$
 $$\left[\gamma_1\in G(e_1,e_2), \gamma_2\in G(e_2.e_3)\right]\Rightarrow \gamma_1 \gamma_2\in G(e_1, e_3).$$
It follows that the union  $\displaystyle G:=\bigcup_{(e_1,e_2)\in E\times E} G(e_1,e_2)$  is  a wide subgroupoid of $\G$. 

For $x\in X$ define  a set $\Gamma(x):=\{\gamma\in\Gamma: (x;\gamma,p(\er(\gamma)))\in \Phi\}$; because  $\Phi(\Gamma\times p(E))=X$ 
this is a family of non empty sets.
Moreover  since $\gamma\in \Gamma(\Phi(\gamma,p(\er(\gamma))))$ the union of all these sets is $\Gamma$. 

We claim that  if $\gamma_1, \gamma_2\in\Gamma(x)$ then $[\gamma_1]=[\gamma_2]$ in $\G/G$. 
Indeed, if $(x;\gamma_1,p(\er(\gamma_1)))\in \Phi$ and $(x;\gamma_2,p(\er(\gamma_2)))\in \Phi$ then $(p(\er(\gamma_2));s(\gamma_2),x)\in \Phi$ and 
$(p(\er(\gamma_2));s(\gamma_2),\gamma_1,p(\er(\gamma_1)))\in \Phi(id\times \Phi)=\Phi(m\times id)$ therefore
$s(\gamma_2) \gamma_1\in G(\er(\gamma_2),\er(\gamma_1)\subset G$ and $[\gamma_1]=[\gamma_2]$.

In this way we can  define $\psi: X\ni x \mapsto [\gamma_x]\in \G/G$, where $\gamma_x$ is any element of $\Gamma(x)$.
Since any $\gamma$ is in $\Gamma(x)$ for some $x$, this mapping is surjective. The equality  $\psi(x)=\psi(y)$ means that
$(x;\gamma_1,p(\er(\gamma_1)))\in \Phi$ and $(y;\gamma_2,p(\er(\gamma_2)))\in \Phi$ and $\gamma_1=\gamma_2 g$ for some $g\in G$, i.e. 
$(p(\er(\gamma_2));g,p(\er(\gamma_1)))\in \Phi$. So $(x;\gamma_2, g, p(\er(\gamma_1)))\in \Phi(m\times id)=\Phi(id \times \Phi)$, 
so exists $z\in X$ such that: $(x;\gamma_2,z)\in \Phi$ and $(z;g, p(\er(\gamma_1)))\in \Phi$ but then $z=p(\er(\gamma_2))$ and $y=x$.
Therefore $\psi$ is bijective.

The last point is the equality $\tilde{m}_G(id\times \psi)=\psi\Phi $. Using definitions of $\tilde{m}_G$ and $\psi$ one gets:
 $$ (A)\,\,:\,(z;\gamma,x)\in \tilde{m}_G(id\times \psi)\iff 
\exists \gamma_1: (x;\gamma_1,p(\er(\gamma_1)))\in \Phi\,,\,\er(\gamma)=\el(\gamma_1)\,,\,z=[\gamma \gamma_1]\in \Gamma/G$$
$$ (B)\,\,:\,(z;\gamma,x)\in\psi \cdot \Phi\iff 
\er(\gamma)=\rho(x)\,,\,
\exists \tilde{\gamma}: (\gamma x;\tilde{\gamma},p(\rho(\gamma x))\in \Phi\,,\,z=[\tilde{\gamma}]\in \Gamma/G,\,(\gamma x:=\Phi(\gamma, x))$$
$(A)\Rightarrow (B)$  by Prop. \ref{action-equiv}:  $\er(\gamma)=\el(\gamma_1)=\rho(x)$, so $(\gamma,x)\in D(\Phi)$ and 
$(\gamma x; \gamma \gamma_1, p(\rho(x)))\in \Phi$ so one can take $\tilde{\gamma}=\gamma \gamma_1$;\\
$(B)\Rightarrow (A)$ since $\er(\gamma)=\rho(x)$ we have  $(\gamma,x)\in D(\Phi)$; then $(x;s(\gamma),\gamma x)\in \Phi$ 
(again by Prop. \ref{action-equiv}) and $(x;s(\gamma)\tilde{\gamma}, p(\rho(x)))\in \Phi$; put $\gamma_1:=s(\gamma)\tilde{\gamma}$ and get $(A)$.\\
\dow

\subsection{ Induced actions}

Let $\Gamma\rra E$ be a groupoid and $\Delta\rra F$ ($F\subset E$)  subgroupoid of $\Gamma$. 
Let $\Phi:\Delta\times X\rel X$ be an action of $\Delta$ on $X$. We assume that its base map $\rho_\Phi: X\rightarrow F$ is surjective 
(otherwise we can restrict the action to the the image of $\rho_\Phi$ which is a union of transistive components of $\Delta$ and 
this is also a subgroupoid of $\Gamma$).
Let $\hat{Y}:=\{(\gamma,x)\in \Gamma\times X: e_R(\gamma)=\rho_\Phi(x)\}$. $\Gamma$ acts on $e_R^{-1}(F)$ 
by the restriction of canonical left action, i.e
we have the action  $l_F:\Gamma\times e_R^{-1}(F)\rel e_R^{-1}(F)$. It is straightforward that 
$l_F\times id_X:\Gamma\times e_R^{-1}(F)\times X\rel e_R^{-1}(F)\times X$ is an action that can be restricted to $\hat{Y}$; this 
restriction will be denoted by $\hat{l}_Y$. On $\hat{Y}$ consider the relation:
$$ (\gamma,x)\simeq (\gamma_1,x_1)\,\,\iff \,\exists \,\delta\in\Delta: \gamma=\gamma_1 s(\delta)\,,\,(x;\delta,x_1)\in\Phi$$
This is an equivalence relation, let $Y$ denotes the set of classes. The following lemma is rather straightforward.
\begin{lem} Let $(\gamma,y_1)\in D(\hat{l}_Y)$. If $[y_1]=[y_2]$ then $(\gamma,y_2)\in D(\hat{l}_Y)$,
$[\hat{l}_Y(\gamma,y_1)]=[\hat{l}_Y(\gamma,y_2)]$ and the definition $l_Y(\gamma,[y]):=[\hat{l}_Y(\gamma,y)]$ gives an action of $\Gamma$ on $Y$.
\end{lem}
\vspace{-1ex}\vspace{-1ex}\dowl
\begin{defi} The action induced from $(\Delta,\Phi)$ is the action introduced in the previous lemma.\end{defi}
\begin{ex1} Let $\Gamma\rightrightarrows E$ be a transitive groupoid, $G$ -- the isotropy group of $e_0\in E$ 
and $G\times Z\ni (g,z)\mapsto g z\in Z$ an action of $G$ on a set $Z$. We can assume that $\Gamma=E\times G\times E$.
Then $\hat{Y}=E\times G\times Z$,  $Y=E\times Z$ and the induced action $\Phi:\G\times (E\times Z)\rel (E\times Z)$ is given by:
$$Gr(\Phi):=\{(e_1, gz;e_1,g,e_2, e_2, z): e_1,e_2\in E, z\in Z, g\in G\}$$
\end{ex1}

In fact this is the most general form of an action of a transitive groupoid on a set:
\begin{prop} Let $\G:=E\times G\times E$ be a transtitive groupoid and $\Phi:\Gamma\times Z\rel Z$ an action. 
There exists $G$ space $\tilde{Z}$ and an equivariant bijection $\Psi:E\times \tilde{Z}\rightarrow Z$ 
(the action of $\Gamma$ on $E\times \tilde{Z}$ is as in the example above.)
\end{prop}
{\em Proof:} Choose $z_0\in Z$ and let $e_0:=\rho(z_0)$, where $\rho:Z\rightarrow E$ is the  base map of the action $\Phi$. 
Let $\tilde{Z}:=\rho^{-1}(e_0)$, define an action of $G$ on $\tilde{Z}$ by: $(g,\tilde{z})\mapsto g\tilde{z}:=\Phi(e_0,g,e_0;\tilde{z})$ and 
$\Psi: E\times \tilde{Z}\rightarrow Z$ by $(z;e,\tilde{z})\in \Psi \iff (z;e,e_G,e_0,\tilde{z})\in \Phi$ ($e_G$ is the neutral element of $G$).
Let $\tilde{\Phi}$ denotes the action induced from $(g,\tilde{z})\mapsto g\tilde{z}$ as in the example above:
$$\tilde{\Phi}:=\{(e_1,g\tilde{z};e_1, g, e_2, e_2,\tilde{z})\}=
\{(e_1,\tilde{z_1}; e_1,g,e_2, e_2,\tilde{z}): (\tilde{z_1}; e_0, g,e_0,\tilde{z})\in \Phi\}$$
We have to prove that $\Psi$ is a bijection and $\Phi (id\times \Psi) =\Psi \tilde{\Phi}$.
Let us compute:
$$(z;e_1,g,e_2,e_3,\tilde{z})\in \Psi\tilde{\Phi}\iff 
e_2=e_3 \,,\, \exists \tilde{z}_3: (z;e_1,e_G,e_0,\tilde{z}_3)\in \Phi\,,\,(\tilde{z}_3;e_0,g,e_0,\tilde{z})\in\Phi\iff$$
$$\iff e_2=e_3\,,\,(z;e_1,e_G,e_0,e_0,g,e_0,\tilde{z})\in \Phi(id\times\Phi)=\Phi(m\times id)\iff e_2=e_3\,,\,(z;e_1,g,e_0,\tilde{z})\in \Phi$$
on the other hand:
$$(z;e_1,g,e_2,e_3,\tilde{z})\in \Phi(id\times \Psi)\iff \exists z_3: (z;e_1, g, e_2,z_3)\in \Phi\,,\,(z_3; e_3,\tilde{z})\in \Psi\iff $$
$$\exists z_3: (z;e_1, g, e_2,z_3)\in \Phi\,,\,(z_3;e_3,e_G,e_0,\tilde{z})\in \Phi\iff $$
$$(z;e_1,g,e_2,e_3,e_G,e_0,\tilde{z})\in \Phi(id \times \Phi)=\Phi(m\times id)\iff e_2=e_3\,,\,(z;e_1,g,e_0,\tilde{z})\in \Phi$$
It remains to show that $\Psi$ is a bijection.
Let  $z:=\Psi(e,\tilde{z})=\Psi(e_1,\tilde{z}_1)$ i.e.
$(z;e,e_G,e_0,\tilde{z})\in \Phi$ and  $(z;e_1,e_G,e_0,\tilde{z}_1)\in \Phi$, therefore  $(\tilde{z}_1; e_0,e_G,e_1,z)\in \Phi$ 
(by prop. \ref{action-equiv}) and $(\tilde{z}_1;e_0,e_G,e_1,e, e_G,e_0,\tilde{z})\in \Phi(id\times \Phi)=\Phi(m\times id)$.
 It follows that $e=e_1$ and then $\tilde{z}=\tilde{z}_1$.\\
For $z\in Z$ there exists $z_1$ such that $(z_1;e_0,e_G,\rho(z),z)\in \Phi$ and then $\rho(z_1)=e_0$ i.e. $z_1\in \tilde{Z}$ and 
$(z; \rho(z), e_G,e_0,z_1)\in \Phi$ i.e $z=\Psi(\rho(z),z_1)$.\dow


\end{document}